\documentclass{article}
\usepackage{amsmath}
\usepackage{amssymb}
\usepackage{graphicx} 
\usepackage{color}
\usepackage{float}
\usepackage{hyperref}
\usepackage{multirow}
\newtheorem{remark}{Remark}

\title{Local interaction simulation approach for the acoustic wave equation with perfectly matched layer}
\author{
Tao Yu
\footnote{School of Mathematics and Physics, Jinggangshan University, Ji'an, Jiangxi 343009, China.
(yutao\_math@163.com)}
\and 
Tailong Jin
\footnote{School of Mathematical Sciences, Institute of Natural Sciences, and MOE-LSC, Shanghai Jiao Tong University, Shanghai 200240, China. (tailongjin0213@163.com)}
\and
Yanfeng Shen 
\footnote{University of Michigan–Shanghai Jiao Tong University Joint Institute, Shanghai Jiao Tong University, Shanghai 200240, China. (yanfeng.shen@sjtu.edu.cn)}
\and 
Lei Zhang
\footnote{School of Mathematical Sciences, Institute of Natural Sciences, and MOE-LSC, Shanghai Jiao Tong University, Shanghai 200240, China.(lzhang2012@sjtu.edu.cn)}
}
\date{}

\begin{document}
\maketitle
\begin{abstract}
Simulation of the acoustic wave equation plays an important role in various applications, including audio engineering, medical imaging, and fluid dynamics. However, the complexity of the propagation medium can pose challenges, such as the infinite computing region and the interface conditions between different media. In this paper, we construct a method for simulating acoustic wave propagation based on the local interaction simulation approach (LISA) and the perfectly matched layer (PML). This method can simulate wave propagation in a finite computing region and overcome the smoothing process at the interface between different media. Numerical examples demonstrate the effectiveness of this approach.
\end{abstract}

{\bf Keywords:}{ acoustic wave equation, local interaction simulation approach (LISA), perfectly matched layer (PML), interface media.}
\section{Introduction}
\label{sec1}

Wave phenomena are prevalent, with acoustic waves being particularly significant due to their various applications. In this paper, we focus on the classical acoustic wave equation.
\begin{equation}\label{acoustic wave equation}
    \frac{\partial^2 u}{\partial t^2} = c^2 \Delta u=  c^2 (\frac{\partial^2 u}{\partial x^2}+\frac{\partial^2 u}{\partial y^2}),
\end{equation}
where \( u \) represents the displacement and \( (c^2)^{-1} = \frac{\rho}{\mu} \) denotes the wave velocity, with \( \rho \) being the density and \( \mu \) the bulk modulus of the propagation medium. When an external force generates a wave, a source term should be added to the right side of the dynamic equilibrium equation.

Wave propagation is an important topic with a rich history across various scientific and engineering disciplines, including seismology, electromagnetism, optics, and acoustics. The Fast Fourier Transform is an effective analytical method for addressing wave transmission, reflection, mode conversion, and nonlinear higher harmonics \cite{Shen-0}. However, analytical approaches can only handle simple structures of infinite size, offering a basic understanding of wave propagation patterns. A more practical way to simulate wave propagation in complex media is through numerical methods. Classical numerical techniques include the Finite Element Method (FEM), Boundary Element Method (BEM), and Finite Difference Method (FDM).

When using classical numerical methods to simulate wave propagation, we encounter two significant challenges. The first issue is that the medium through which the wave propagates may be complex and inhomogeneous. Accurately simulating the wave propagation process requires the correct implementation of interface conditions; if these conditions are not well-defined, the numerical method may become unstable. Among existing numerical methods, the Finite Difference Method (FDM) is straightforward to derive, implement, and parallelize. However, its main drawback is accuracy, particularly when using smooth parameters at discontinuous interfaces \cite{LISA3}. To address this, Delsanto et al. proposed a powerful modeling technique called the Local Interaction Simulation Approach (LISA), based on explicit finite difference methods and a sharp interface model, for simulating wave propagation in isotropic, heterogeneous materials \cite{LISA1, LISA2, LISA3}. This approach enforces the continuity of displacements and stresses at the interface, enabling precise treatment of perfect interfaces between different materials, making LISA more accurate than pure FDM in scenarios involving discontinuous material properties. The accuracy of LISA was numerically justified in \cite{Sundararaman2}. Additionally, LISA was later implemented in parallel using graphics processing units (GPUs) with Compute Unified Device Architecture (CUDA), demonstrating that GPU-based LISA outperformed FEM and reduced computation time from hours to minutes \cite{Paćko}. In recent years, numerous scholars have conducted extensive research related to this topic \cite{Lee1, Lee2, Sundararaman1, Cesnik-1, Cesnik-2, Radecki, Shen-1, Shen-2, Shen-3}.

Another challenge in numerical simulations of acoustic waves is that they typically occur in infinite or semi-infinite domains. However, our primary focus is often on the local propagation characteristics of these waves. To address this, we should employ non-reflective boundary techniques (NRBs). Existing NRBs can be categorized into two types, depending on the need for an extended region: absorbing boundary conditions (ABCs) and absorbing boundary layers (ABLs). ABCs impose conditions directly on the boundaries of a structure, while ABLs require an extended region connected to the boundaries to dampen incident waves and minimize reflections. ABCs were introduced by B. Engquist and A. Majda for the acoustic wave equation \cite{Engquist}. They are particularly effective for absorbing perpendicular incident waves; however, for waves traveling obliquely to the boundary, higher-order ABCs are necessary to achieve acceptable accuracy \cite{Hagstrom}. While absorbing boundary conditions have been utilized in finite element methods \cite{Xue}, the issue of boundary reflection remains unresolved. Reflection becomes pronounced when waves travel obliquely to the boundary, and deriving the control equations for higher-order ABCs in FEM can be challenging. For ABLs, two common approaches are the perfectly matched layer (PML) \cite{Berenger} and the absorbing layer using increasing damping (ALID) \cite{Israeli}. Both methods aim to dampen waves in an extended absorbing region, differing in whether the impedance of the simulated region matches that of the extended region. The PML was initially introduced by B\'{e}renger for the absorption of 2-D electromagnetic waves \cite{Berenger}, and Chew and Weedon demonstrated its effectiveness using stretched coordinates \cite{Chew}. Variants of PML, such as complex frequency shifted perfectly matched layer (CFS-PML) \cite{Kuzuoglu} and auxiliary differential equation perfectly matched layer (ADE-PML) \cite{Ramadan}, have also been proposed. Generally, PML is derived from first-order equations, with the second-order wave equation transformed into first-order form. Michael provides a general PML format for wave equations in \cite{PML-wave}, while a more rigorous derivation of second-order PML is presented in this paper.

The LISA method was initially developed for the elastic wave equation with a focus on the engineering community. In this paper, we present a concise mathematical derivation of the LISA algorithm in the context of the acoustic wave equation. We also derive the first-order and second-order perfectly matched layer (PML) conditions within the PML framework. Following this, we integrate the LISA method with the PML conditions to simulate acoustic wave propagation. Additionally, we implement the method on a GPU and demonstrate its effectiveness in simulating the acoustic wave propagation.

The outline of the paper is as follows: In Section \ref{sec2}, we provide a brief description of the Local Interaction Simulation Approach for the acoustic wave equation. Section \ref{sec3} presents the construction of the first-order and second-order perfectly matched layers. Finally, Section \ref{sec4} includes several numerical experiments.

\section{Local Interaction Simulation Approach}
\label{sec2}

In this section, we introduce the general strategy for the LISA approach, in particular, for the acoustic wave equation. 

To derive the numerical formula, we first need to discretize the space and time domains. Assume the spatial domain is divided into a grid of \( N_x \times N_y \) rectangular cells, each of size \( h_x \times h_y \). The unit time step for time discretization is represented by \( \tau \). If the medium is homogeneous, we can employ the classical finite difference method to simulate wave propagation, such as the second-order explicit central difference scheme.
\begin{eqnarray}
    u_{i,j}^{t+1} = 2 u_{i,j}^{t} - u_{i,j}^{t-1} +  \tau^2c^2 ( \!\!\!\!&&\!\!\!\! \frac{u_{i+1,j}^{t}-2u_{i,j}^{t}+u_{i-1,j}^{t}}{h_x^2}\nonumber\\ \!\!\!\!&+&\!\!\!\!\frac{u_{i,j+1}^{t}-2u_{i,j}^{t}+u_{i,j-1}^{t}}{h_y^2}),
\end{eqnarray}
where $u_{i,j}^{t}$ denotes the numerical value of displacement $u$ at the grid point $(i,j)$ and time step $t$. For simplicity, we assume equal grid lengths in both spatial directions, defined as $h_x = h_y \triangleq h$.

\begin{figure}[H]
    \centering
    \includegraphics[width=0.9\linewidth]{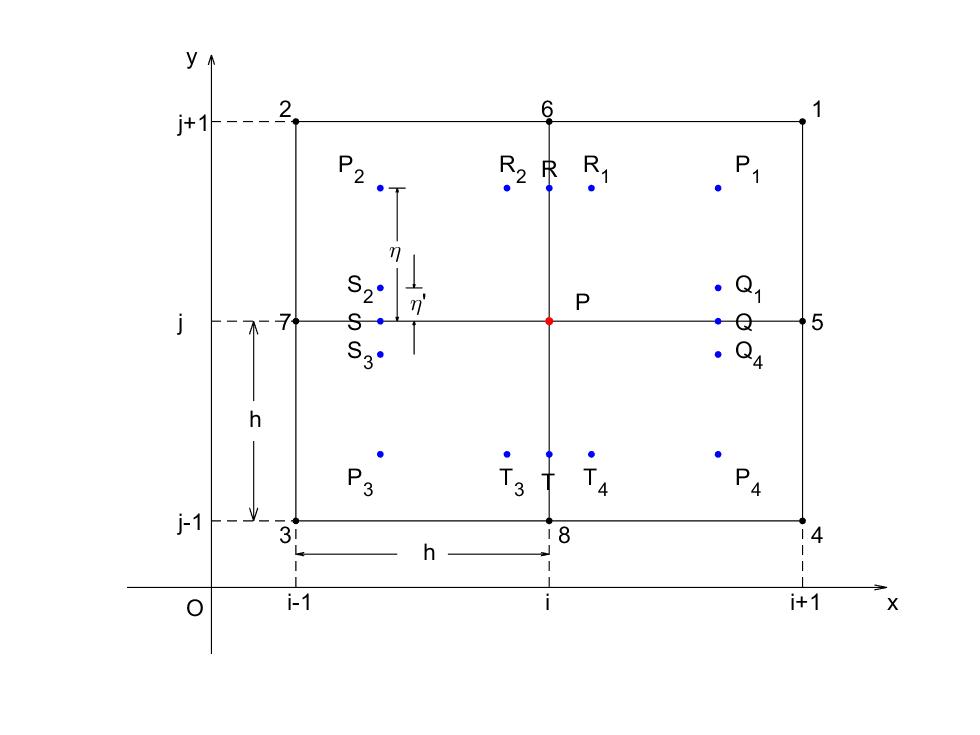}
    \caption{Graphical illustration of the cross point $P$ and other auxiliary points used in the derivation of LISA.} \label{crosspoint}  
\end{figure}

If the medium is inhomogeneous, several assumptions must be made to apply the Local Interaction Simulation Approach to the acoustic wave equation. After spatial discretization, we assume that the physical properties of the medium are constant within each cell but may vary from cell to cell. We then need to derive iteration formulas for the cross point $P$ (see Figure \ref{crosspoint}), which is located at the intersection of two orthogonal interfaces separating four media with different physical properties $\rho_k, \mu_k$ (for $k=1,\cdots,4$). The Local Interaction Simulation Approach requires us to impose two continuity conditions at the cross point $P$.

\paragraph{Continuity condition I:} The first condition is the continuity of the second time derivative of \( u \) at cross point $P$. To apply this condition, we construct four auxiliary points \( P_k \) (for \( k=1,\cdots,4 \)) near the cross point \( P \). The distance between each \( P_k \) and the grid is denoted as \( \eta \) (where \( 0 < \eta \ll h \)). As \( \eta \rightarrow 0 \), we assume that the second time derivative of the displacement computed at each \( P_k \) converges to a common value \( \Omega \),
\begin{equation}\label{condition1}
   \frac{\partial^2 u}{\partial t^2}\bigg|_{P_k} \approx \Omega.
\end{equation}
$\Omega$ can be viewed as the second time derivative of the displacement $u$ at the cross point $P$. Let $u_i\;(i=1,2,\cdots,8)$ represent the value of $u$ at the point labeled $i$ (see Figure \ref{crosspoint}) and applying Taylor's formula to the displacement $u$ in the point labeled 5 at the point $P_1$, we get
\begin{eqnarray}
    u_5 \!\!\!\!&\approx&\!\!\!\! u_{P_1} + u_{P_1,x}(h-\eta)+ u_{P_1,y}\eta + \frac{1}{2} u_{P_1,xx}(h-\eta)^2 +u_{P_1,xy}\eta(h-\eta)  + \frac{1}{2} u_{P_1,yy}\eta^2\nonumber\\ \!\!\!\!&\approx&\!\!\!\! u + u_{P_1,x} h + \frac{1}{2} u_{P_1,xx} h^2,
\end{eqnarray}
as $\eta \rightarrow 0$.

Thus
\begin{equation}\label{u_xx}
    u_{P_1,xx} \approx -\frac{2}{h} u_{P_1,x} + \frac{2}{h^2} (u_5-u).
\end{equation}
Similarly, 
\begin{equation}\label{u_yy}
    u_{P_1,yy} \approx -\frac{2}{h} u_{P_1,y} + \frac{2}{h^2} (u_6-u).
\end{equation}

Substituting equations (\ref{condition1}) - (\ref{u_yy}) into the wave equation (\ref{acoustic wave equation}) at point $P_1$, we get
\begin{eqnarray}\label{Omega}
 \Omega \!\!\!\!&\approx&\!\!\!\! c_1^2 \left( -\frac{2}{h} u_{P_1,x} + \frac{2}{h^2} (u_5-u) -\frac{2}{h} u_{P_1,y} + \frac{2}{h^2} (u_6-u) \right) \nonumber\\
\!\!\!\!&\approx&\!\!\!\! -\frac{2 c_1^2}{h} u_{P_1, x} - \frac{2 c_1^2}{h} u_{P_1, y} + \frac{2 c_1^2}{h^2} (u_5-u) + \frac{2 c_1^2}{h^2} (u_6-u).
\end{eqnarray}
Multiply both sides of the above equation (\ref{Omega}) by $h^2$, it becomes
\begin{equation}\label{condition11}
    h^2 \Omega \approx -2 c_1^2 h u_{P_1,x} - 2 c_1^2 h u_{P_1,y} +2 c_1^2 (u_5-u) + 2 c_1^2 (u_6-u).
\end{equation}
For the other three points, we can similarly get
\begin{eqnarray}
    h^2 \Omega \!\!\!\!&\approx&\!\!\!\! \;\;\;2 c_2^2 h u_{P_2,x} - 2 c_2^2 h u_{P_2,y} +2 c_2^2 (u_7-u) + 2 c_2^2 (u_6-u),\\\label{condition12}
    h^2 \Omega \!\!\!\!&\approx&\!\!\!\! \;\;\;2 c_3^2 h u_{P_3,x} + 2 c_3^2 h u_{P_3,y} +2 c_3^2 (u_7-u) + 2 c_3^2 (u_8-u),\\\label{condition13}
    h^2 \Omega \!\!\!\!&\approx&\!\!\!\! - 2 c_4^2 h u_{P_4,x} + 2 c_4^2 h u_{P_4,y} +2 c_4^2 (u_5-u) + 2 c_4^2 (u_8-u).\label{condition14}
\end{eqnarray}

\paragraph{Continuity condition II:} the second condition is the continuity of the flux $f=\frac{1}{\rho}\nabla u\cdot \vec{n}$ at four interfaces connecting $P$, which is
\begin{equation}\label{condition2}
    f_{R_2} = f_{R_1},\; f_{S_3} = f_{S_2},\; f_{T_4} = f_{T_3},\; f_{Q_1} = f_{Q_4},
\end{equation}
where $\vec{n}$ is the outer normal vector at the interface and $R_1, R_2, S_2, S_3, T_3, T_4, Q_4, Q_1$ are the other eight required auxiliary points near the interface, see Figure \ref{crosspoint}. The distance of these points from the nearest interface is denoted as $\eta' (0<\eta'\ll \eta)$. For example, the first one of the continuity conditions (\ref{condition2}) is equivalent to
\begin{eqnarray}\label{condition21}
     f_{R_2} = f_{R_1} \!\!\!\!&\Leftrightarrow&\!\!\!\! \left(\frac{1}{\rho}\nabla u\cdot \vec{n}\right)_{R_2}=\left(\frac{1}{\rho}\nabla u\cdot \vec{n}\right)_{R_1}\Leftrightarrow
    \frac{1}{\rho_2}u_{R_2,x}=\frac{1}{\rho_1}u_{R_1,x}\nonumber\\
    \!\!\!\!&\Rightarrow&\!\!\!\!
    \frac{1}{\rho_2}u_{P_2,x}\approx\frac{1}{\rho_1}u_{P_1,x}\Leftrightarrow
    \frac{1}{\rho_2}u_{P_2,x}-\frac{1}{\rho_1}u_{P_1,x}\approx 0.
\end{eqnarray}
similarly, we can obtain the other conditions, 
\begin{eqnarray}
   \frac{1}{\rho_3} u_{P_3,x} - \frac{1}{\rho_4} u_{P_4,x}\approx 0,\\\label{condition22}
   \frac{1}{\rho_4} u_{P_4,y} - \frac{1}{\rho_1} u_{P_1,y}\approx 0,\\\label{condition23}
   \frac{1}{\rho_3} u_{P_3,y} - \frac{1}{\rho_2} u_{P_2,y}\approx 0.\label{condition24}
\end{eqnarray}
Summing the above equations (\ref{condition21}) - (\ref{condition24}) together, we derive
\begin{eqnarray}
   \!\!\!\!&-&\!\!\!\!\frac{1}{\rho_1}u_{P_1,x}-\frac{1}{\rho_1} u_{P_1,y} + \frac{1}{\rho_2}u_{P_2,x} - \frac{1}{\rho_2} u_{P_2,y} \nonumber\\
   \!\!\!\!&+&\!\!\!\!\frac{1}{\rho_3} u_{P_3,x} + \frac{1}{\rho_3} u_{P_3,y} -\frac{1}{\rho_4} u_{P_4,x} + \frac{1}{\rho_4} u_{P_4,y} \approx 0.\label{condition15}
\end{eqnarray}
On the other hand, after dividing by \( 2\rho c^2 = 2\mu \), equations (\ref{condition11}) through (\ref{condition14}) can be re-formulated as follows:
\begin{eqnarray}
    \frac{h^2 \Omega}{2\mu_1} \!\!\!\!&\approx&\!\!\!\! - \frac{h}{\rho_1} u_{P_1,x} - \frac{h}{\rho_1} u_{P_1,y} +\frac{1}{\rho_1} (u_5-u) + \frac{1}{\rho_1} (u_6-u),\\
    \frac{h^2 \Omega}{2\mu_2} \!\!\!\!&\approx&\!\!\!\! \;\;\;\frac{h}{\rho_2} u_{P_2,x} - \frac{h}{\rho_2} u_{P_2,y} +\frac{1}{\rho_2} (u_7-u) + \frac{1}{\rho_2} (u_6-u),\\
    \frac{h^2 \Omega}{2\mu_3} \!\!\!\!&\approx&\!\!\!\! \;\;\;\frac{h}{\rho_3} u_{P_3,x} + \frac{h}{\rho_3} u_{P_3,y} +\frac{1}{\rho_3} (u_7-u) + \frac{1}{\rho_3} (u_8-u),\\
    \frac{h^2 \Omega}{2\mu_4} \!\!\!\!&\approx&\!\!\!\! - \frac{h}{\rho_4} u_{P_4,x} + \frac{h}{\rho_4} u_{P_4,y} +\frac{1}{\rho_4} (u_5-u) + \frac{1}{\rho_4} (u_8-u).
\end{eqnarray}

Using \ref{condition15}, the summation of the four equations above leads to,
\begin{eqnarray}
    (\frac{1}{\mu_1}+\frac{1}{\mu_2}+\frac{1}{\mu_3}+\frac{1}{\mu_4}) \frac{h^2}{2} \Omega \approx \!\!\!\!&&\!\!\!\! h (-\frac{1}{\rho_1}u_{P_1,x}-\frac{1}{\rho_1} u_{P_1,y} + \frac{1}{\rho_2}u_{P_2,x} - \frac{1}{\rho_2} u_{P_2,y}\nonumber\\ 
    \!\!\!\!&+&\!\!\!\! \frac{1}{\rho_3} u_{P_3,x} + \frac{1}{\rho_3} u_{P_3,y} -\frac{1}{\rho_4} u_{P_4,x} + \frac{1}{\rho_4} u_{P_4,y})\nonumber\\
    \!\!\!\!&+&\!\!\!\! (\frac{1}{\rho_1}+\frac{1}{\rho_4})u_5 + (\frac{1}{\rho_1}+\frac{1}{\rho_2})u_6+(\frac{1}{\rho_2}+\frac{1}{\rho_3})u_7\nonumber\\
    \!\!\!\!&+&\!\!\!\! (\frac{1}{\rho_3}+\frac{1}{\rho_4})u_8-2(\frac{1}{\rho_1}+\frac{1}{\rho_2}+\frac{1}{\rho_3}+\frac{1}{\rho_4})u\nonumber\\
    \approx \!\!\!\!&&\!\!\!\! (\frac{1}{\rho_1}+\frac{1}{\rho_4})u_5 + (\frac{1}{\rho_1}+\frac{1}{\rho_2})u_6+(\frac{1}{\rho_2}+\frac{1}{\rho_3})u_7\nonumber\\
    \!\!\!\!&+&\!\!\!\! (\frac{1}{\rho_3}+\frac{1}{\rho_4})u_8-2(\frac{1}{\rho_1}+\frac{1}{\rho_2}+\frac{1}{\rho_3}+\frac{1}{\rho_4})u,\label{sumation equation}    
\end{eqnarray}

Denote by $\alpha=\frac{1}{4}\sum_{k=1}^{4}\frac{1}{\mu_k},\; \rho_5=[\frac{1}{2}(\frac{1}{\rho_1}+\frac{1}{\rho_4})]^{-1},\; \rho_6=[\frac{1}{2}(\frac{1}{\rho_1}+\frac{1}{\rho_2})]^{-1},\; \rho_7=[\frac{1}{2}(\frac{1}{\rho_2}+\frac{1}{\rho_3})]^{-1},\; \rho_8=[\frac{1}{2}(\frac{1}{\rho_3}+\frac{1}{\rho_4})]^{-1},\; \beta=(\sum_{k=1}^{4}\frac{1}{\rho_k})^{-1}, $ the previous equation (\ref{sumation equation}) is equivalent to
\begin{equation}
    \alpha h^2 \Omega \approx -\frac{1}{\beta} u+\sum_{k=5}^{8}\frac{1}{\rho_k}u_k.
\end{equation}

Using the second order scheme in time,
\begin{equation}
    \Omega=\frac{\partial^2 u}{\partial t^2}\bigg|_{P}\approx \frac{u^{t+1}-2u+u^{t-1}}{\tau^2},
\end{equation}
we obtain the final LISA algorithm
\begin{equation}\label{LISA}
    u^{t+1}=2u-u^{t-1}+\frac{\tau^2}{h^2\alpha}\left(-\frac{1}{\beta} u+\sum_{k=5}^{8}\frac{1}{\rho_k}u_k\right).
\end{equation}
\begin{remark}(Homogeneous Case)
    If the medium is homogeneous, namely, $\alpha=\mu^{-1}, \beta=\rho/4$ and $\rho_k=\rho\;(k=5,\cdots,8)$. The LISA algorithm (\ref{LISA}) becomes
    \begin{equation}
        u^{t+1}=2u-u^{t-1}+\frac{\tau^2}{h^2}\frac{\mu}{\rho}\left(u_5+u_7+u_6+u_8-4u\right),
    \end{equation}
    which reduces to the classical second order explicit finite difference scheme for acoustic wave equation.
\end{remark}
\begin{remark}(Multilayer Case)
    If the medium has only a vertical interface, for example, assume $\mu_1=\mu_4 \triangleq \mu_r,\; \mu_2=\mu_3\triangleq\mu_l$ and $\rho_1=\rho_4\triangleq\rho_r,\; \rho_2=\rho_3\triangleq\rho_l$. The LISA algorithm (\ref{LISA}) reduces to
    \begin{eqnarray}
        u^{t+1}=2u-u^{t-1}+\frac{\tau^2}{h^2}\frac{\mu_l \mu_r}{\mu_l+ \mu_r}\frac{1}{\rho_l \rho_r}(-4(\rho_l+\rho_r)u\!\!\!\!&+&\!\!\!\!2\rho_l u_5+2\rho_r u_7\nonumber\\
        \!\!\!\!&+&\!\!\!\!(\rho_l+\rho_r)(u_6+u_8)).
    \end{eqnarray}
\end{remark}

\section{Perfectly Matched Layer}
\label{sec3}
In this section, we first briefly introduce the perfectly matched layer (PML) method for the acoustic wave equation, as described in \cite{Abarbanel, Berenger, Chew, Duru, PML-wave}. We then present a rigorous derivation of the second-order PML using a complex-valued coordinate stretching technique. Consequently, we derive the PML condition for the LISA scheme to simulate wave propagation in a finite domain.

\subsection{PML formulation}

The acoustic wave equation (\ref{acoustic wave equation}) can be rewritten in divergence form as follows:
\begin{equation}\label{divergence form equation}
    \frac{\partial^2 u}{\partial t^2} = \mu \nabla \cdot \left(\frac{1}{\rho} \nabla u\right).
\end{equation}
By introducing the intermediate varibles \( v = (v_x, v_y) \), the above equation (\ref{divergence form equation}) is equivalent to a system of first-order equations.
\begin{eqnarray}
\label{first order equations}
\left\{
\begin{array}{ll}\displaystyle
\frac{\partial u}{\partial t} =\mu \nabla\cdot v,\\\displaystyle
\frac{\partial v}{\partial t} = \frac{1}{\rho}\nabla u.
\end{array}
\right.
\end{eqnarray}

\paragraph{PML equation for $v$:} We first focus on the second equation of the system (\ref{first order equations}). It can be rewritten in scalar form as follows,
\begin{equation}
    \frac{\partial v_x}{\partial t}=\frac{1}{\rho}\frac{\partial u}{\partial x},\;\;\frac{\partial v_y}{\partial t}=\frac{1}{\rho}\frac{\partial u}{\partial y}.
\end{equation}
Fourier transform of the above equations can be expressed as,
\begin{equation}
    -i\omega v_x^*=\frac{1}{\rho}\frac{\partial u^*}{\partial x},\;\;-i\omega v_y^*=\frac{1}{\rho}\frac{\partial u^*}{\partial y},
\end{equation}
where \( \omega \) denotes the circular frequency, and \( u^* \) and \( v^* = (v_x^*, v_y^*) \) represent the Fourier transform of displacement and intermediate variables. Following the PML theory, the original coordinates \( x \) and \( y \) are transformed into complex coordinates \( \tilde{x} \) and \( \tilde{y} \). Consequently, the above equations become,
\begin{equation}\label{complex coordinates equation}
    -i\omega v^*_x=\frac{1}{\rho}\frac{\partial u^*}{\partial \tilde{x}},\;\;-i\omega v^*_y=\frac{1}{\rho}\frac{\partial u^*}{\partial \tilde{y}},
\end{equation}
where the complex coordinate stretching is defined by
\begin{equation}
\tilde{x}=x+\frac{i}{\omega}\int_{0}^{x}\lambda_x(s)dx,\;\;\tilde{y}=y+\frac{i}{\omega}\int_{0}^{y}\lambda_y(s)dx,
\label{eqn:complexcoordinatestretching}
\end{equation}
and $\lambda_x(x)$ and $\lambda_y(y)$ are the damping functions in $x$ and $y$ directions. The partial derivatives for $\tilde{x}$, $\tilde{y}$ and the partial derivatives for $x$, $y$ have relations
\begin{eqnarray}
\frac{\partial}{\partial \tilde{x}}=(1+\frac{i\lambda_x}{\omega})^{-1}\frac{\partial}{\partial x},\;\;\frac{\partial}{\partial \tilde{y}}=(1+\frac{i\lambda_y}{\omega})^{-1}\frac{\partial}{\partial y}.
\end{eqnarray}
Plug the above equations into equation (\ref{complex coordinates equation}), it becomes 
\begin{equation}
    -i\omega v^*_x=(1+\frac{i\lambda_x}{\omega})^{-1}\frac{1}{\rho}\frac{\partial u^*}{\partial x},\;\;-i\omega v^*_y=(1+\frac{i\lambda_y}{\omega})^{-1}\frac{1}{\rho}\frac{\partial u^*}{\partial y}.
\end{equation}
Multiplying both sides of the above equations by $1+\frac{i\lambda_x}{\omega}$ and  $1+\frac{i\lambda_y}{\omega}$, separately, we can derive 
\begin{equation}
    -i\omega v^*_x=-\lambda_x v^*_x+\frac{1}{\rho}\frac{\partial u^*}{\partial x},\;\;-i\omega v^*_y=-\lambda_y v^*_y+\frac{1}{\rho}\frac{\partial u^*}{\partial y}.
\end{equation}
By inverse Fourier transform, we obtain the time domain equation, 
\begin{equation}
    \frac{\partial v_x}{\partial t} = -\lambda_x v_x+\frac{1}{\rho}\frac{\partial u}{\partial x},\;\; \frac{\partial v_y}{\partial t}=-\lambda_y v_y+\frac{1}{\rho}\frac{\partial u}{\partial y}.
\end{equation}
Denote by $\lambda=(\lambda_x,\lambda_y)$, the vector form reads
\begin{equation}\label{PML1}
    \frac{\partial v}{\partial t} =  \frac{1}{\rho} \nabla u-\lambda \odot v ,
\end{equation}
where the vector operator $\odot$ is defined by $\lambda\odot v=(\lambda_x v_x, \lambda_y v_y)$ for vectors $\lambda=(\lambda_x, \lambda_y)$ and $v=(v_x, v_y)$. 

\paragraph{PML equation for $u$} Next, we focus on the first equation of the system (\ref{first order equations}). By applying the Fourier transform to shift everything into the frequency domain, we obtain,
\begin{equation}
    -i\omega u^* = \mu \nabla\cdot v^*.
\end{equation}
By the complex coordinate stretching \ref{eqn:complexcoordinatestretching}, the above equation can be transformed to 
\begin{equation}
    -i\omega u^* = \mu [(1+\frac{i\lambda_x}{\omega})^{-1}\frac{\partial v^*_x}{\partial x}+(1+\frac{i\lambda_y}{\omega})^{-1}\frac{\partial v^*_y}{\partial y}].
\end{equation}
Multiplying $(1+\frac{i\lambda_x}{\omega})(1+\frac{i\lambda_y}{\omega})$ to both sides of the above equation, we can obtain
\begin{equation}
    -i\omega u^*(1+\frac{i\lambda_x}{\omega})(1+\frac{i\lambda_y}{\omega}) = \mu [(1+\frac{i\lambda_y}{\omega})\frac{\partial v^*_x}{\partial x}+(1+\frac{i\lambda_x}{\omega})\frac{\partial v^*_y}{\partial y}].
\end{equation}
After simple calculations, this equation is equivalent to
\begin{equation}
    -i\omega u^* + (\lambda_x+\lambda_y)u^* = \mu \nabla \cdot v^*+\frac{i}{\omega}(\mu\lambda_y\frac{\partial v^*_x}{\partial x}+\mu\lambda_x\frac{\partial v^*_y}{\partial y}-\lambda_x\lambda_y u^*),
\end{equation}

By introducing an auxiliary variable $\psi$, such that its Fourier transform $\psi^*$ satisfies
\begin{equation}
    \psi^*=\frac{i}{\omega}(\mu\lambda_y\frac{\partial v^*_x}{\partial x}+\mu\lambda_x\frac{\partial v^*_y}{\partial y}-\lambda_x\lambda_y u^*),
\end{equation}
and multiplying both sides by $-i\omega$, we can get 
\begin{equation}
    -i\omega\psi^*=\mu\lambda_y\frac{\partial v^*_x}{\partial x}+\mu\lambda_x\frac{\partial v^*_y}{\partial y}-\lambda_x\lambda_y u^*.
\end{equation}
Bringing it back to the time domain, we have
\begin{equation}\label{PML12}
    \frac{\partial \psi}{\partial t}=\mu\lambda_y\frac{\partial v_x}{\partial x}+\mu\lambda_x\frac{\partial v_y}{\partial y}-\lambda_x\lambda_y u.
\end{equation}
Then, the PML equation for $u$ becomes 
\begin{equation}\label{PML13}
\frac{\partial u}{\partial t}=\mu \nabla \cdot v-(\lambda_x+\lambda_y)u+\psi.
\end{equation}

\paragraph{PML system} Putting everything together, we have the following PML system of equations,
\begin{eqnarray}\label{1orderpml}
\left\{
\begin{array}{ll}\displaystyle
\frac{\partial v}{\partial t} = \frac{1}{\rho} \nabla u -\lambda \odot v,\\\displaystyle
\frac{\partial \psi}{\partial t}=\mu\lambda_y\frac{\partial v_x}{\partial x}+\mu\lambda_x\frac{\partial v_y}{\partial y}-\lambda_x\lambda_y u,\\\displaystyle
\frac{\partial u}{\partial t}=\mu \nabla \cdot v-(\lambda_x+\lambda_y)u+\psi.
\end{array}
\right.
\end{eqnarray}
It has one more equation than the original wave equation (\ref{first order equations}). Since this is a first-order equation of displacement, it is so-called \emph{first-order method}.  One can use the explicit finite difference method to solve the above three differential equations. At each time step, the finite difference scheme for the first equation of (\ref{1orderpml}) is
\begin{eqnarray}\label{FDM1-1order}
   (v_{x})_{i,j}^{t+1} \!\!\!\!&=&\!\!\!\! (v_{x})_{i,j}^{t}+\tau\left(-\lambda_x(x_i)(v_x)_{i,j}^{t}+\frac{u_{i+1,j}^t-u_{i,j}^t}{\rho h}\right),\\ \label{FDM2-1order}
    (v_{y})_{i,j}^{t+1} \!\!\!\!&=&\!\!\!\! (v_{y})_{i,j}^{t}+\tau\left(-\lambda_y(y_j)(v_y)_{i,j}^{t}+\frac{u_{i,j+1}^t-u_{i,j}^t}{\rho h}\right).
\end{eqnarray}
The finite difference scheme for the other two equations of (\ref{1orderpml}) are
\begin{eqnarray}
   \psi_{i,j}^{t+1} = \psi_{i,j}^{t}+\tau( \!\!\!\!& &\!\!\!\! \mu\lambda_y(y_j)\frac{(v_x)_{i+1,j}^{t+1}-(v_x)_{i,j}^{t+1}}{h}\nonumber\\
   \!\!\!\!&+&\!\!\!\! \mu\lambda_x(x_i)\frac{(v_y)_{i,j+1}^{t+1}-(v_y)_{i,j}^{t+1}}{h} -\lambda_x(x_i)\lambda_y(y_j)u_{i,j}^t)
\end{eqnarray}
and
\begin{eqnarray}
   u_{i,j}^{t+1} = u_{i,j}^{t}+\tau( \!\!\!\!& &\!\!\!\!  \mu\frac{(v_x)_{i+1,j}^{t+1}-(v_x)_{i,j}^{t+1}}{h}+\mu\frac{(v_y)_{i,j+1}^{t+1}-(v_y)_{i,j}^{t+1}}{h} \nonumber\\
   \!\!\!\!&-&\!\!\!\!(\lambda_x(x_i)+\lambda_y(y_j))u_{i,j}^{t} + \phi_{i,j}^{t+1}).
\end{eqnarray}

\subsection{Second Order PML Method}

The so-called second-order method can be derived in the same way. Equation (\ref{PML1}) for $v$, derived from the second equation of the system (\ref{first order equations}), is still needed. However, the treatment of the first equation of (\ref{first order equations}) is different. Taking the partial derivative with respect to time $t$, the first equation of (\ref{first order equations}) becomes
\begin{equation}\label{deformation-equation}
    \frac{\partial^2 u}{\partial t^2} =\mu \nabla\cdot \frac{\partial v}{\partial t}.
\end{equation}

Following the PML theory and utilizing equation (\ref{PML1}), we obtain
\begin{equation}
    -\omega^2 u^* = \mu [(1+\frac{i\lambda_x}{\omega})^{-1}\frac{\partial}{\partial x}(-\lambda_x v^*_x+\frac{1}{\rho}\frac{\partial u^*}{\partial x})+(1+\frac{i\lambda_y}{\omega})^{-1}\frac{\partial}{\partial y}(-\lambda_y v^*_y+\frac{1}{\rho}\frac{\partial u^*}{\partial y})],
\end{equation}
by applying Fourier transform and complex coordinate stretching. Multiplying both sides of the above equation by $(1+\frac{i\lambda_x}{\omega})(1+\frac{i\lambda_y}{\omega})$, we have
\begin{eqnarray}
    -\omega^2 (1+\frac{i\lambda_x}{\omega})(1+\frac{i\lambda_y}{\omega}) u^* = \mu [ \!\!\!\!& &\!\!\!\! (1+\frac{i\lambda_y}{\omega})\frac{\partial}{\partial x}(-\lambda_x v^*_x+\frac{1}{\rho}\frac{\partial u^*}{\partial x})\nonumber\\
     \!\!\!\!&+&\!\!\!\!(1+\frac{i\lambda_x}{\omega})\frac{\partial}{\partial y}(-\lambda_y v_y+\frac{1}{\rho}\frac{\partial u^*}{\partial y})].
     \label{eqn:secondorderPML}
\end{eqnarray}

After applying the inverse Fourier transform (IFT), the left-hand side can be expressed as:
\begin{eqnarray}
    -\omega^2 (1+\frac{i\lambda_x}{\omega})(1+\frac{i\lambda_y}{\omega}) u^*
    \!\!\!\!&=&\!\!\!\!-\omega^2u^*- i\omega(\lambda_x+\lambda_y)u^* + \lambda_x \lambda_y u^* \nonumber\\
\!\!\!\!&\xrightarrow{IFT}&\!\!\!\! \frac{\partial^2 u}{\partial t^2}+(\lambda_x+\lambda_y)\frac{\partial u}{\partial t}+\lambda_x \lambda_y u
\end{eqnarray}

The right hand side of \ref{eqn:secondorderPML} is more involved, and can be calculated as,
\begin{eqnarray}
    \!\!\!\!& &\!\!\!\!\mu [(1+\frac{i\lambda_y}{\omega})\frac{\partial}{\partial x}(-\lambda_x v^*_x+\frac{1}{\rho}\frac{\partial u^*}{\partial x})+(1+\frac{i\lambda_x}{\omega})\frac{\partial}{\partial y}(-\lambda_y v^*_y+\frac{1}{\rho}\frac{\partial u^*}{\partial y})]\nonumber\\
     \!\!\!\!&=&\!\!\!\!\frac{\mu}{\rho}\nabla^2 u^*-\mu\nabla\cdot(\lambda \odot v^*)\nonumber\\
     \!\!\!\!&&\!\!\!\!+\frac{i}{\omega}\mu[\frac{\lambda_y}{\rho}\frac{\partial^2 u^*}{\partial x^2}+\frac{\lambda_x}{\rho}\frac{\partial^2 u^*}{\partial y^2}-\lambda_y\frac{\partial}{\partial x}(\lambda_x v^*_x)-\lambda_x\frac{\partial}{\partial y}(\lambda_y v^*_y)],
\end{eqnarray}
We introduce an auxiliary variation $\psi$, defined by using its Fourier transform $\psi^*$
\begin{equation}
    \psi^*=\frac{i}{\omega}\mu[\frac{\lambda_y}{\rho}\frac{\partial^2 u^*}{\partial x^2}+\frac{\lambda_x}{\rho}\frac{\partial^2 u^*}{\partial y^2}-\lambda_y\frac{\partial}{\partial x}(\lambda_x v^*_x)-\lambda_x\frac{\partial}{\partial y}(\lambda_y v^*_y)].
\end{equation}
Multiplying both sides by $-i\omega$, we obtain
\begin{equation}
    -i\omega\psi^* =\mu[\frac{\lambda_y}{\rho}\frac{\partial^2 u^*}{\partial x^2}+\frac{\lambda_x}{\rho}\frac{\partial^2 u^*}{\partial y^2}-\lambda_y\frac{\partial}{\partial x}(\lambda_x v^*_x)-\lambda_x\frac{\partial}{\partial y}(\lambda_y v^*_y)].
\end{equation}
Returning to the time domain, we arrive at the following equation,
\begin{equation}\label{PML2}
    \frac{\partial \psi}{\partial t}=\mu[\frac{\lambda_y}{\rho}\frac{\partial^2 u}{\partial x^2}+\frac{\lambda_x}{\rho}\frac{\partial^2 u}{\partial y^2}-\lambda_y\frac{\partial}{\partial x}(\lambda_x v_x)-\lambda_x\frac{\partial}{\partial y}(\lambda_y v_y)].
\end{equation}
at the same time, equation (\ref{deformation-equation}) transforms into 
\begin{equation}\label{PML3}
\frac{\partial^2 u}{\partial t^2}+(\lambda_x+\lambda_y)\frac{\partial u}{\partial t}+\lambda_x \lambda_y u=\frac{\mu}{\rho}\nabla^2 u-\mu\nabla\cdot(\lambda \odot v)+\psi.
\end{equation}
Therefore, we solve the following second order system in PML media
\begin{eqnarray}\label{2orderpml}
\left\{
\begin{array}{ll}\displaystyle
\frac{\partial v}{\partial t} = \frac{1}{\rho} \nabla u-\lambda \odot v,\\\displaystyle
\frac{\partial \psi}{\partial t}=\mu[\frac{\lambda_y}{\rho}\frac{\partial^2 u}{\partial x^2}+\frac{\lambda_x}{\rho}\frac{\partial^2 u}{\partial y^2}-\lambda_y\frac{\partial}{\partial x}(\lambda_x v_x)-\lambda_x\frac{\partial}{\partial y}(\lambda_y v_y)],\\\displaystyle
\frac{\partial^2 u}{\partial t^2}+(\lambda_x+\lambda_y)\frac{\partial u}{\partial t}+\lambda_x \lambda_y u=\frac{\mu}{\rho}\nabla^2 u-\mu\nabla\cdot(\lambda \odot v)+\psi.
\end{array}
\right.
\end{eqnarray}

We can similarly construct an explicit difference scheme to solve the equations above. At each time step, the finite difference scheme for the first equation of (\ref{2orderpml}) can be chosen to be the same as that of the first-order method, represented by equations (\ref{FDM1-1order}) and (\ref{FDM2-1order}). For the second equation of (\ref{2orderpml}), we can apply the Euler forward difference scheme,
\begin{eqnarray}
   \psi_{i,j}^{t+1} = \psi_{i,j}^{t}+\tau\mu(\!\!\!\!&&\!\!\!\!\frac{\lambda_y(y_j)}{\rho}\cdot\frac{u^{t}_{i+1,j}-2u^{t}_{i,j}+u^{t}_{i-1,j}}{h^2}\nonumber\\
   \!\!\!\!&+&\!\!\!\! 
   \frac{\lambda_x(x_i)}{\rho}\cdot\frac{u^{t}_{i,j+1}-2u^{t}_{i,j}+u^{t}_{i,j-1}}{h^2}\nonumber\\   \!\!\!\!&-&\!\!\!\!\lambda_y(y_j)\frac{\lambda_x(x_{i+1})(v_x)_{i+1,j}^{t+1}-\lambda_x(x_{i-1})(v_x)_{i-1,j}^{t+1}}{2h}\nonumber\\   \!\!\!\!&-&\!\!\!\!\lambda_x(x_i)\frac{\lambda_y(y_{j+1})(v_y)_{i,j+1}^{t+1}-\lambda_y(y_{j-1})(v_y)_{i,j-1}^{t+1}}{2h}).
\end{eqnarray}

By applying central differences to the time derivatives in the third equation of (\ref{2orderpml}), we get,
\begin{eqnarray}   
\!\!\!\!&&\!\!\!\!\frac{u^{t+1}_{i,j}-2u^{t}_{i,j}+u^{t-1}_{i,j}}{\tau^2}+(\lambda_x(x_i)+\lambda_y(y_j))\frac{u^{t+1}_{i,j}-u^{t-1}_{i,j}}{2\tau} +\lambda_x(x_i)\lambda_y(y_j)u^{t}_{i,j}
\nonumber\\ \!\!\!\!&=&\!\!\!\!-\mu(\frac{\lambda_x(x_{i+1})(v_x)^{t+1}_{i+1,j}-\lambda_x(x_{i-1})(v_x)^{t+1}_{i-1,j}}{2h}+\frac{\lambda_y(y_{j+1})(v_y)^{t+1}_{i,j+1}-\lambda_y(y_{j-1})(v_y)^{t+1}_{i,j-1}}{2h}) \nonumber\\
\!\!\!\!&+&\!\!\!\!\frac{\mu}{\rho}(\frac{u^{t}_{i+1,j}-2u^{t}_{i,j}+u^{t}_{i-1,j}}{h^2}+\frac{u^{t}_{i,j+1} -2u^{t}_{i,j}+u^{t}_{i,j-1}}{h^2})+\psi_{i,j}^{t+1}.
\end{eqnarray}
which is essentially an explicit scheme
\begin{eqnarray}   
\!\!\!\!&(&\!\!\!\!1+\frac{\tau}{2}(\lambda_x(x_i)+\lambda_y(y_j)))u^{t+1}_{i,j}=(2-\lambda_x(x_i)\lambda_y(y_j))u^{t}_{i,j}-(1-\frac{\tau}{2}(\lambda_x(x_i)+\lambda_y(y_j)))u^{t-1}_{i,j}
\nonumber\\ \!\!\!\!&-&\!\!\!\!\frac{\mu \tau^2}{2h}(\lambda_x(x_{i+1})(v_x)^{t+1}_{i+1,j}-\lambda_x(x_{i-1})(v_x)^{t+1}_{i-1,j}+\lambda_y(y_{j+1})(v_y)^{t+1}_{i,j+1}-\lambda_y(y_{j-1})(v_y)^{t+1}_{i,j-1}) \nonumber\\
\!\!\!\!&&\!\!\!\!+\frac{\mu \tau^2}{\rho h^2}(u^{t}_{i+1,j}-2u^{t}_{i,j}+u^{t}_{i-1,j}+u^{t}_{i,j+1} -2u^{t}_{i,j}+u^{t}_{i,j-1})+ \tau^2\psi_{i,j}^{t+1}.
\end{eqnarray}

All that remains is to define the decay functions \( \lambda_x(x) \) and \( \lambda_y(y) \). A straightforward approach is to use a simple polynomial, 
\begin{equation}
    \lambda_x(x) = \lambda_{max}(x-x_0)^m,
\end{equation}
where $x_0$ is the point at which the PML is activated, and $\lambda_{\text{max}}$ is the maximum value required for proper attenuation, defined by,
\begin{equation}
    \lambda_{max}=-\frac{c\log(R)(m+1)}{(x_{max}-x_{0})^{m+1}}.
\end{equation}
Depending on accuracy requirements, $R$ is typically set to $10^{-4}$ and $m$ is set to 4. The width of the PML is given by $|x_{\text{max}} - x_{0}|$. The decay function $\lambda_y(y)$ can be defined similarly.

\section{Numerical Experiments}
\label{sec4}

In this section, we will present several numerical examples to evaluate the performance of the LISA method.

\subsubsection*{Example 1: Constant Coefficient and Dirichlet BC} Consider the acoustic wave equation (\ref{acoustic wave equation}) with constant coefficients. Assume the density \( \rho = 1 \) and the bulk modulus \( \mu = 1 \). The computational domain \( [0,1]^2 \) is discretized using a uniform grid. To observe the convergence rate of the LISA scheme (\ref{LISA}), the spatial steps are chosen as \( h = \frac{1}{64}, \frac{1}{128}, \frac{1}{256}, \frac{1}{512} \). The temporal step is \( \Delta t = 0.1 h \), which is sufficiently small to ensure stability. Meanwhile, the temporal error can be considered negligible compared to the spatial error. For simplicity, the final simulation time is set to \( T = 1 \). Using Dirichlet boundary conditions, the numerical solution is compared with an analytic solution given by 
\begin{equation}\label{Exact-solution}
    U(x,y,t) = \cos(10\pi x + 1) \cos(10\pi y + 2) \cos(10\pi \sqrt{2} t + 3),
\end{equation}
as used in \cite{FDM_wave}. The relative error is defined as \( \frac{||u_h - U_h||}{||U_h||} \), where \( u_h \) is the numerical solution and \( U_h \) is the exact solution projected onto the spatial grid with step size \( h \). The convergence rate is computed by 
\begin{equation}
 \mathrm{order} = \frac{\log\left(\frac{||u_{2h} - U_{2h}||/||U_{2h}||}{||u_h - U_h||/||U_h||}\right)}{\log(2)}.
\end{equation}

Table \ref{table1} presents the relative errors and convergence rates under the discrete \( L^{2} \) and \( L^{\infty} \) norms. The definitions of these norms are as follows,
\[
||u||_{L^2} = \sqrt{h_x h_y \sum_{j=1}^{N_y+1} \sum_{i=1}^{N_x+1} |u_{i,j}|^2},
\]
and
\[
||u||_{L^{\infty}} = \max_{1 \leq i \leq N_x+1, 1 \leq j \leq N_y+1} |u_{i,j}|.
\]

It is evident that the LISA scheme with constant coefficients behaves like a second-order explicit central difference scheme. One advantage of using the LISA formula is the ability to leverage GPUs to accelerate runtime. Table \ref{table2} presents the program runtime with and without CUDA. The main program is written in Python, using Numba for the CUDA implementation. For the constant coefficients case, the effect of using a GPU to accelerate runtime is not immediately apparent. However, as the grid scale increases, the acceleration effect becomes more pronounced.

\begin{table}[H]
\centering
\begin{tabular}{|c|c|c|c|c|c|c|}
\hline
 N    &$L^2$ error & order  &$L^{\infty}$ error & order\\\hline
 64    &1.8286e-01  & -      &3.7630e-01         & -    \\\hline
 128  &4.9521e-02  & 1.8846 &1.0243e-01         & 1.8772\\\hline
 256  &1.2705e-02  & 1.9626 &2.5980e-02         & 1.9792\\\hline
 512  &3.1993e-03  & 1.9896 &6.5277e-03         & 1.9927\\\hline
\end{tabular}
\caption{Relative error results and convergence orders for acoustic wave equation with constant coefficients}\label{table1}
\end{table}

\begin{table}[H]
\centering
\begin{tabular}{|c|c|c|}
\hline
 N    &Runtime without CUDA (s)  & Runtime with CUDA (s)  \\\hline
 64   & 0.370209      &  2.947894      \\\hline
 128  & 1.551154      &  2.137946      \\\hline
 256  & 8.543154      &  5.420279      \\\hline
 512  &68.042135      & 17.377175      \\\hline
\end{tabular}
\caption{Runtime with and without CUDA in Example 2}\label{table2}
\end{table}

\subsubsection*{Example 2: Non-constant Coefficient and Dirichlet BC} 

Consider a similar case to Example 1, but with non-constant coefficients in the acoustic wave equation (\ref{acoustic wave equation}). The density and the bulk modulus are chosen as \( \rho = 1 \) and \( \mu = 1 + \frac{1}{2} \cos(2\pi x) \cos(2\pi y) \). Consequently, the coefficient \( c^2(x,y) \) is given by 
\begin{equation}
c^2(x,y) = 1 + \frac{1}{2} \cos(2\pi x) \cos(2\pi y).
\end{equation}
For simplicity, we consider the same exact solution as in (\ref{Exact-solution}) for this example. The acoustic wave equation then becomes
\begin{equation}
    \frac{\partial^2 u(x,y)}{\partial t^2} =  c^2(x,y) \left( \frac{\partial^2 u(x,y)}{\partial x^2} + \frac{\partial^2 u(x,y)}{\partial y^2} \right) + f(x,y,t),
\end{equation}
where the source function is 
\begin{equation}
    f(x,y,t) = (10\pi)^2 \cos(2\pi x) \cos(2\pi y) \cos(10\pi x + 1) \cos(10\pi y + 2) \cos(10\pi \sqrt{2} t + 3).
\end{equation}

The relative errors and convergence orders are shown in Table \ref{table3}. It is evident that the convergence orders are reduced to order one. This occurs because the parameter values in each cell are replaced by the parameter values at the center of each cell. Table \ref{table4} presents the runtime with and without CUDA. For the non-constant coefficients case, the acceleration effect of using a GPU is significant.

\begin{table}[H]
\centering
\begin{tabular}{|c|c|c|c|c|c|c|}
\hline
 N    &$L^2$ error & order  &$L^{\infty}$ error & order\\\hline
 64   &1.8752e-01  & -      &3.3960e-01         & -    \\\hline
 128  &5.6847e-02  & 1.7219 &9.3494e-02         & 1.8609\\\hline
 256  &2.2024e-02  & 1.3680 &4.0457e-02         & 1.2085\\\hline
 512  &9.9918e-03  & 1.1403 &1.8607e-02         & 1.1205\\\hline
\end{tabular}
\caption{Relative error results and convergence orders for acoustic wave equation with non-constant coefficients}\label{table3}
\end{table}

\begin{table}[H]
\centering
\begin{tabular}{|c|c|c|}
\hline
 N    &Runtime without CUDA (s)  & Runtime with CUDA (s)  \\\hline
 64   &  30.143023      & 4.147659      \\\hline
 128  & 226.393971      & 3.464801      \\\hline
 256  & 1053.614328     & 11.251251      \\\hline
 512  & 8442.370489     & 37.067896      \\\hline
\end{tabular}
\caption{Runtime with and without CUDA in Example 2}\label{table4}
\end{table}

\subsubsection*{Example 3: Interface Case} 
Consider the multilayer case in which the medium has only a vertical interface. The computational region is \( [-50, 50]^2 \), and the interface is located at \( x = 25 \). The final time is chosen to be \( T = 40 \), allowing the wave to pass through the interface without reaching the boundary. The density and bulk modulus on the left and right sides of the interface are \( \rho_l = 1, \; \mu_l = 1 \) and \( \rho_r = 2, \; \mu_r = 2 \), respectively.

To generate an acoustic wave, we use the following source function 
\begin{equation}\label{force}
    f(x,y,t) = \frac{1}{2\pi} e^{-\frac{x^2+y^2}{2}} \cos t,
\end{equation}
for the first \( [\pi/\Delta t] \) time steps, as used in \cite{PML-wave}. Since there is no analytical expression for the exact solution, the numerical solution derived by the LISA scheme on a uniform grid with spatial step \( h = \frac{1}{4096} \) and temporal step \( \Delta t = 0.1 h \) is considered as the exact solution. The LISA scheme is applied with spatial steps \( h = \frac{1}{128}, \frac{1}{256}, \frac{1}{512}, \frac{1}{1024} \). The relative errors and convergence rates are shown in Table \ref{table5}. From the results, we can see that the LISA method handles the multilayer case very effectively.

\begin{table}[H]
\centering
\begin{tabular}{|c|c|c|c|c|c|c|}
\hline
 N    &$L^2$ error & order  &$L^{\infty}$ error & order\\\hline
 128  & 6.5292e-01  &    -    & 6.7801e-01 &    -   \\\hline 
 256  & 2.2878e-01  & 1.5129  & 3.1559e-01 & 1.1032 \\\hline
 512  & 5.9132e-02  & 1.9520  & 7.7086e-02 & 2.0335 \\\hline
 1024 & 1.5046e-02  & 1.9746  & 2.0021e-02 & 1.9450 \\\hline
\end{tabular}
\caption{Relative error results and convergence orders for the multilayer case}\label{table5}
\end{table}

\subsubsection*{Example 4: Outgoing Wave with Constant Coefficient and PML} 
The acoustic wave equation is computed in the domain \( [-50, 50]^2 \). We begin by discretizing the domain using a uniform grid with \( h = \frac{1}{256} \) for convenience. To ensure stability, we select the time step \( \Delta t = \frac{h}{4c} \), where \( c \) is the wave velocity. 

First, we observe the propagation of an acoustic wave through a medium with density \( \rho = 1 \) and bulk modulus \( \mu = 1 \). To generate the acoustic wave, we apply the same source function as in (\ref{force}) during the first \( [\pi/\Delta t] \) time steps. In the perfectly matched layer (PML), we choose 15 spatial steps. 

\begin{figure}[H]
	\centering
	\begin{minipage}{0.28\linewidth}
		\centering
		\includegraphics[width=1\linewidth]{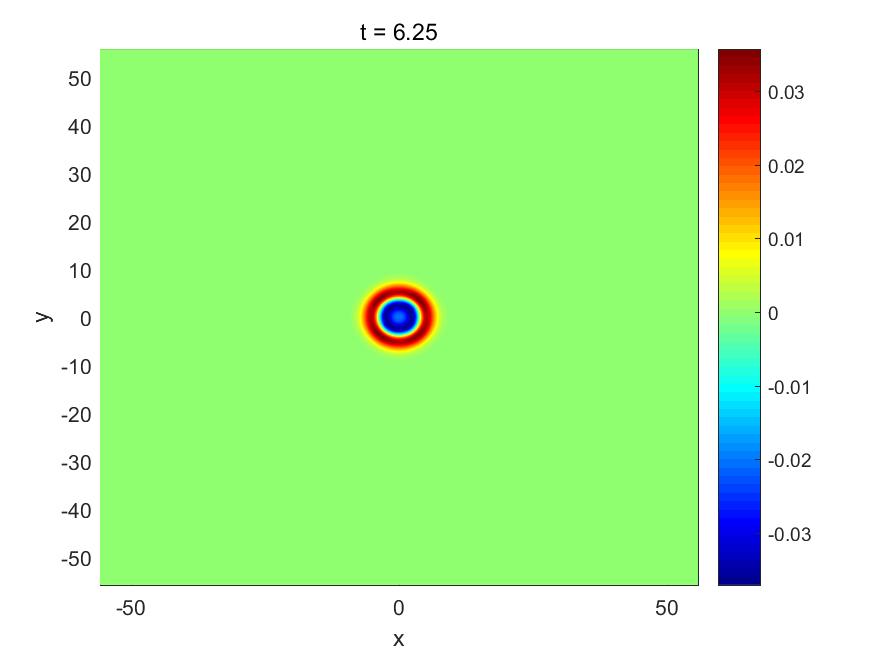}
	\end{minipage}
	\begin{minipage}{0.28\linewidth}
		\centering
		\includegraphics[width=1\linewidth]{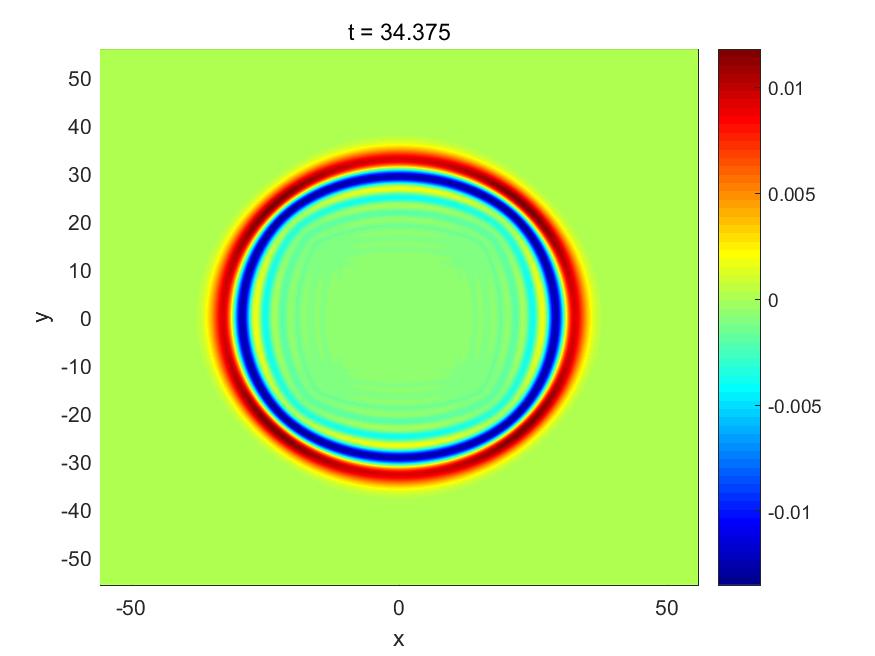} 
	\end{minipage}
    \begin{minipage}{0.28\linewidth}
		\centering
		\includegraphics[width=1\linewidth]{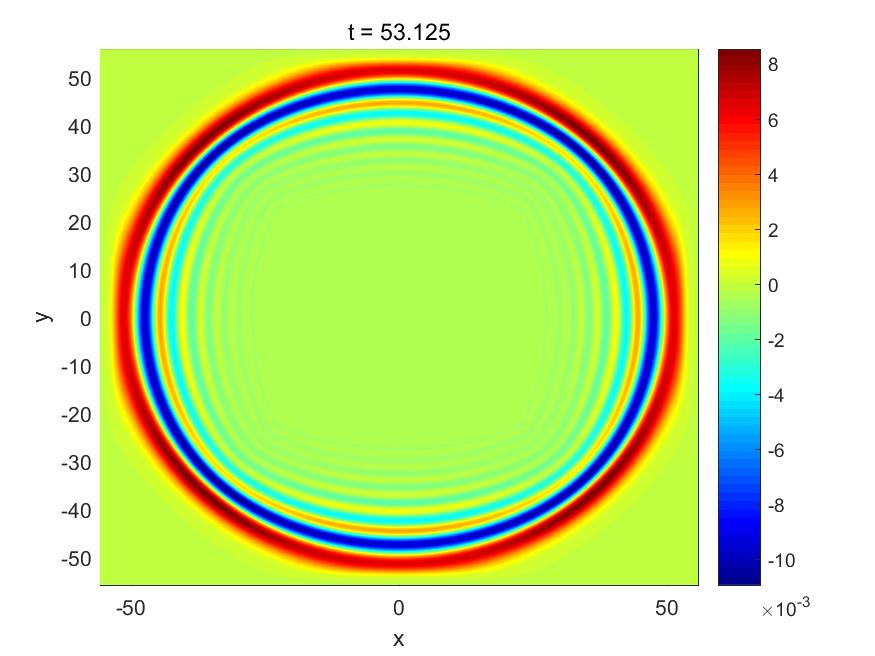}
	\end{minipage}
    \begin{minipage}{0.28\linewidth}
		\centering
		\includegraphics[width=1\linewidth]{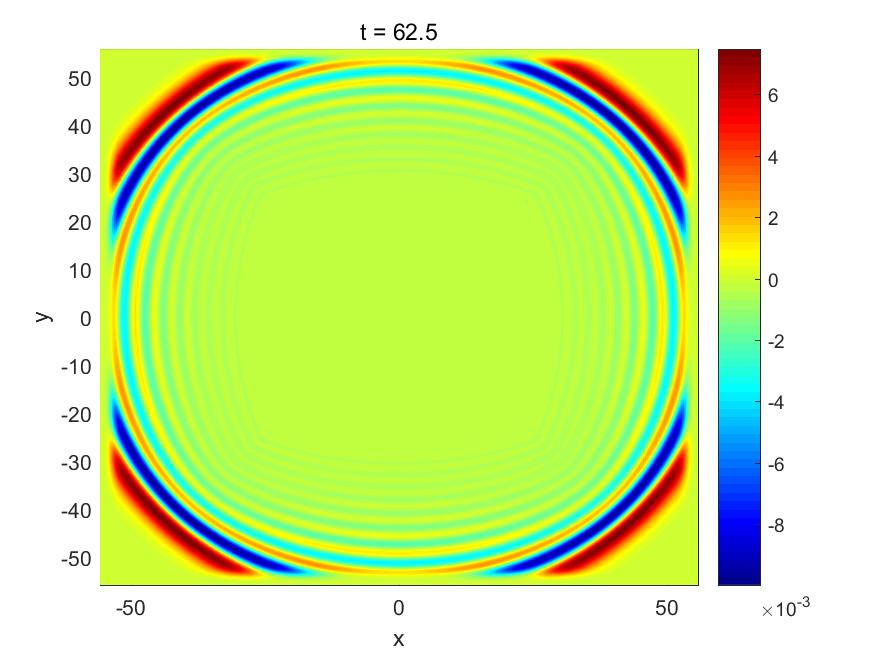}
	\end{minipage} 
    \begin{minipage}{0.28\linewidth}
		\centering
		\includegraphics[width=1\linewidth]{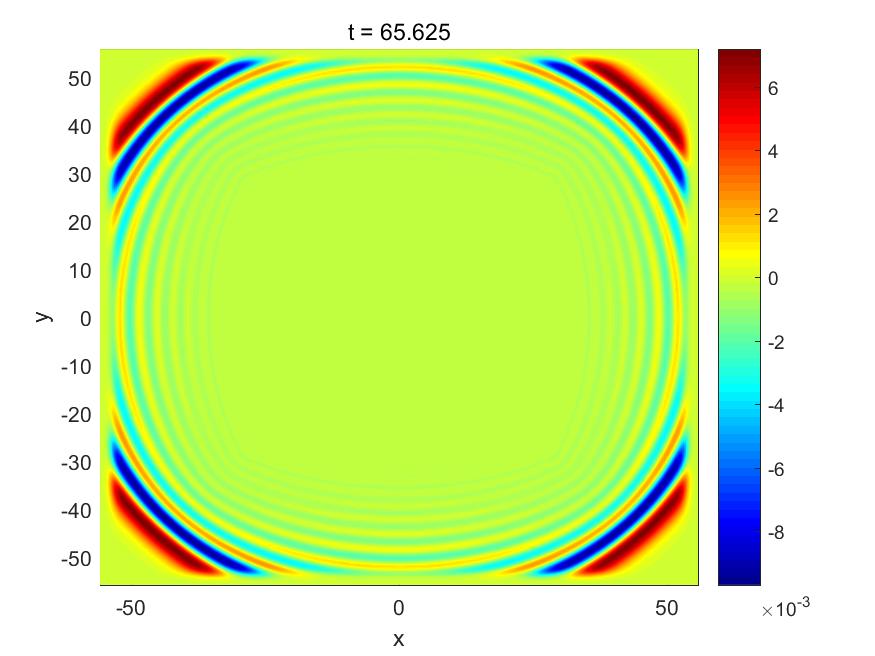}
	\end{minipage}
    \begin{minipage}{0.28\linewidth}
		\centering
		\includegraphics[width=1\linewidth]{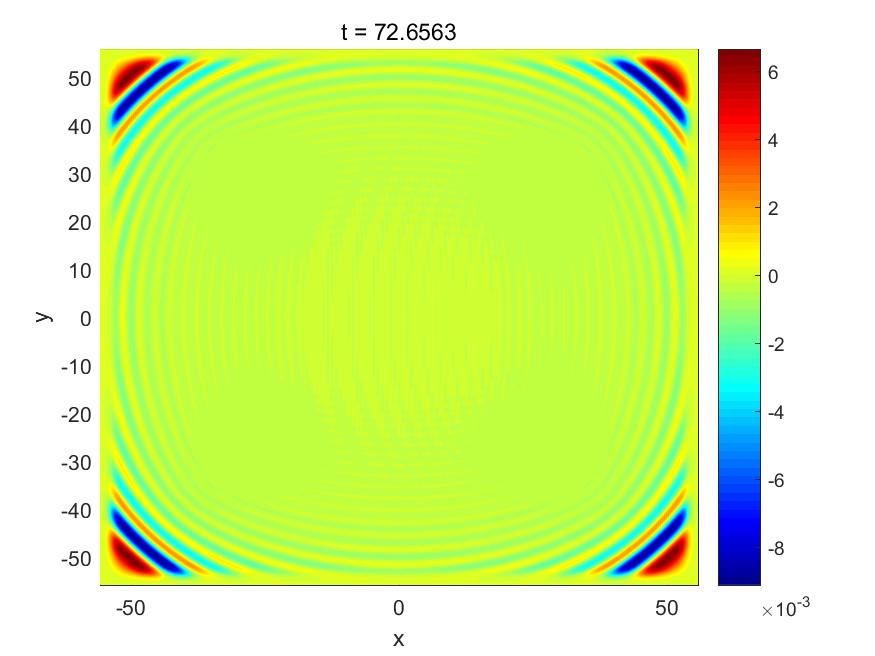}
	\end{minipage} 
 \caption{The snapshot of wave propagation without interface}
 \label{The snapshot of wave propagation with PML}
\end{figure}

Snapshots of the acoustic wave propagation process are shown in Figure \ref{The snapshot of wave propagation with PML}, where the first-order PML (\ref{1orderpml}) is used. As observed in the figures, the wave is effectively absorbed in the perfectly matched layer. Using the second-order PML (\ref{2orderpml}) yields similar results. 

For a more intuitive look at the difference between first-order PML and second-order PML, we compare the relative errors when using these two PMLs, separately. The numerical solution with spatial grid size $N = 4096$ is viewed as an exact solution. The final time is chosen to be $T = 60$. We derive numerical solutions with spatial steps $h = \frac{1}{128}, \frac{1}{256}, \frac{1}{512}, \frac{1}{1024}$, separately. Relative errors are shown in Table \ref{table6}. It can be easily seen that these errors do not vary greatly. For the acoustic wave equation, the first-order PML is sufficient for absorption.

For a more intuitive comparison between first-order and second-order perfectly matched layers (PML), we analyze the relative errors associated with each. We consider the numerical solution with a spatial grid size of \( N = 4096 \) as the exact solution, with the final time set to \( T = 60 \). We derive numerical solutions using spatial steps of \( h = \frac{1}{128}, \frac{1}{256}, \frac{1}{512}, \frac{1}{1024} \). The relative errors are presented in Table \ref{table6}, showing that these errors do not vary significantly. For the acoustic wave equation, the first-order PML proves sufficient for effective absorption.

\begin{table}[H]
\centering
\begin{tabular}{|c|c|c|c|c|}
\hline
\multicolumn{1}{|c|}{\multirow{2}{*}{N}}&\multicolumn{2}{|c|}{$L^2$ error}&\multicolumn{2}{|c|}{$L^{\infty}$ error}\\\cline{2-5}
      & first-order PML & second-order PML & first-order PML & second-order PML\\\hline
 128  & 8.8643e-01  & 8.8786e-01 & 6.5625e-01 & 6.5852e-01 \\\hline 
 256  & 2.6080e-01  & 2.6407e-01 & 2.5158e-01 & 2.6381e-01 \\\hline
 512  & 6.6115e-02  & 7.5837e-02 & 5.8758e-02 & 7.2577e-02 \\\hline
 1024 & 5.9743e-02  & 6.3876e-02 & 3.4909e-02 & 4.3118e-02 \\\hline
\end{tabular}
\caption{Relative error results with first-order PML and second-order PML}\label{table6}
\end{table}

\textbf{Example 5: Outgoing Wave with Interface Coefficients and PML} 
We observe the propagation of an outgoing wave through interface media. For convenience, we assume that the interface between the two media is perpendicular to the X-axis at \( x = 25 \). The densities of the left and right media are \( \rho_l = 1 \) and \( \rho_r = 2 \), respectively. The bulk moduli are \( \mu_l = 1 \) and \( \mu_r = 2 \), respectively. 

\begin{figure}[H]
	\centering
	\begin{minipage}{0.28\linewidth}
		\centering
		\includegraphics[width=1\linewidth]{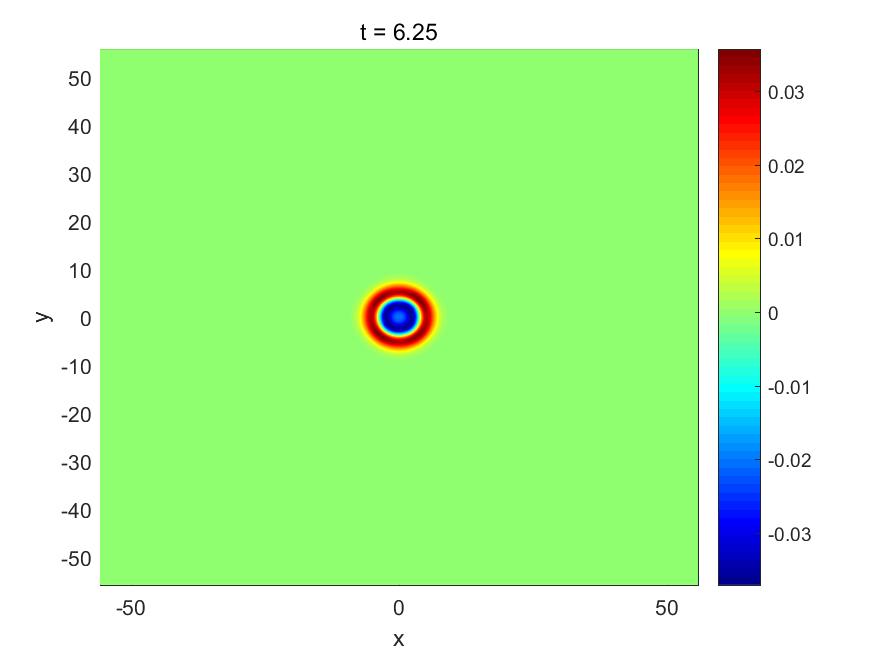}
	\end{minipage}
	\begin{minipage}{0.28\linewidth}
		\centering
		\includegraphics[width=1\linewidth]{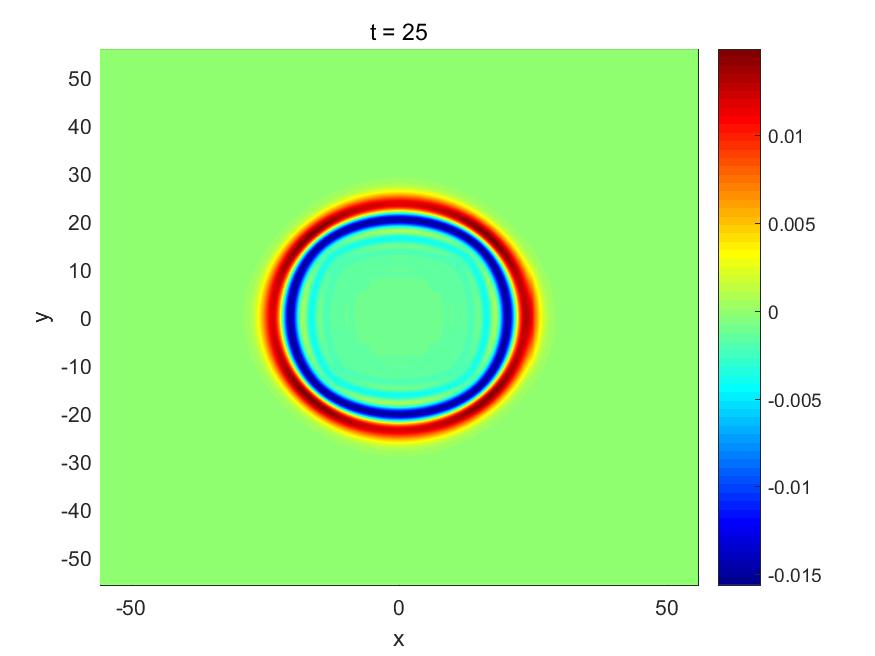} 
	\end{minipage}
    \begin{minipage}{0.28\linewidth}
		\centering
		\includegraphics[width=1\linewidth]{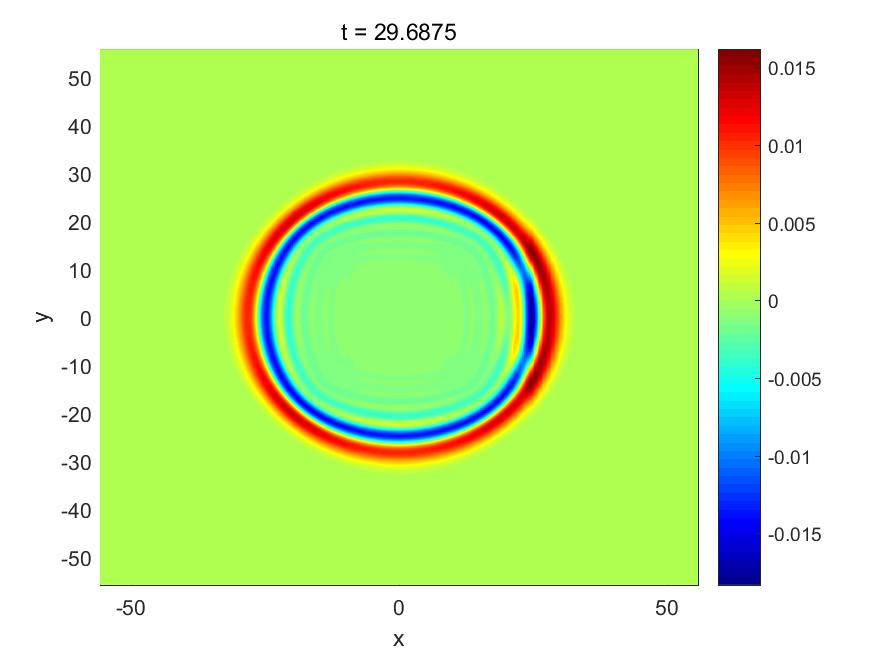}
	\end{minipage}
    \begin{minipage}{0.28\linewidth}
		\centering
		\includegraphics[width=1\linewidth]{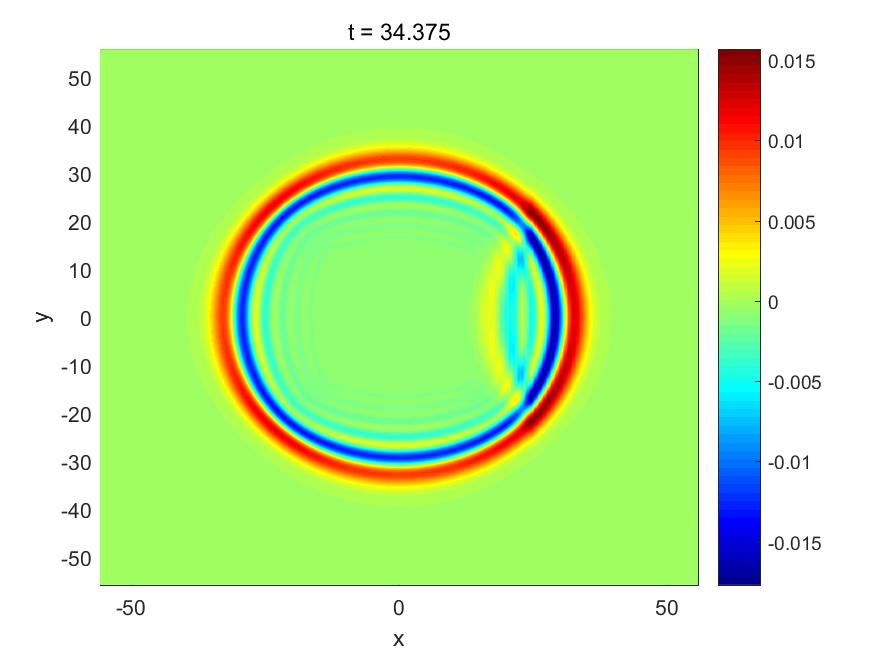}
	\end{minipage} 
    \begin{minipage}{0.28\linewidth}
		\centering
		\includegraphics[width=1\linewidth]{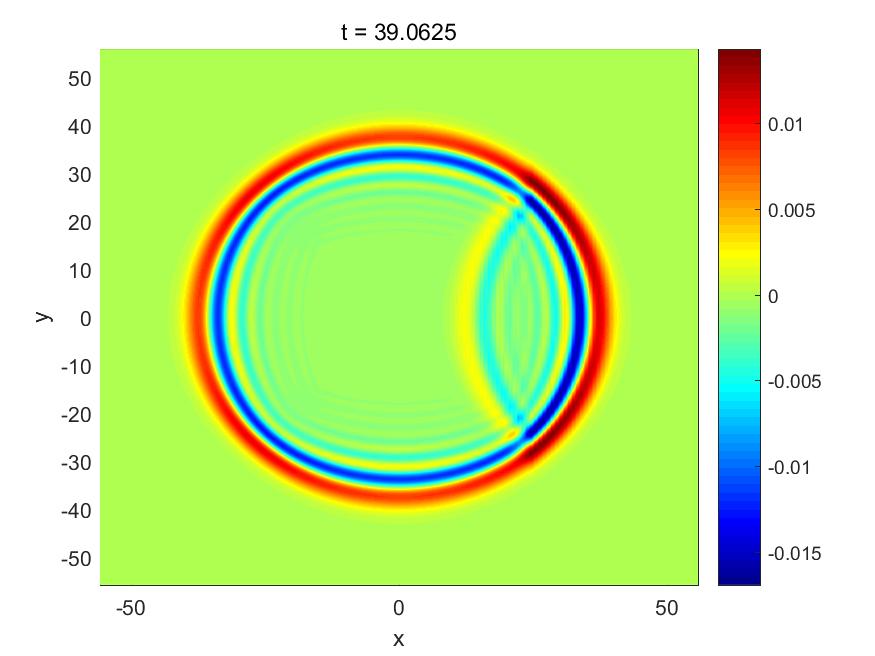}
	\end{minipage}
    \begin{minipage}{0.28\linewidth}
		\centering
		\includegraphics[width=1\linewidth]{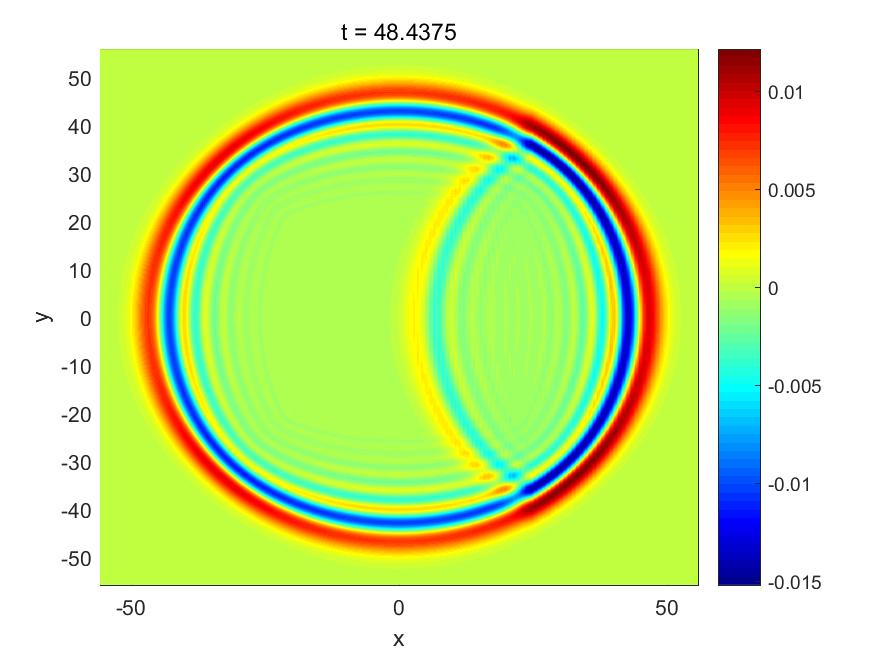}
	\end{minipage} 
     \begin{minipage}{0.28\linewidth}
		\centering
		\includegraphics[width=1\linewidth]{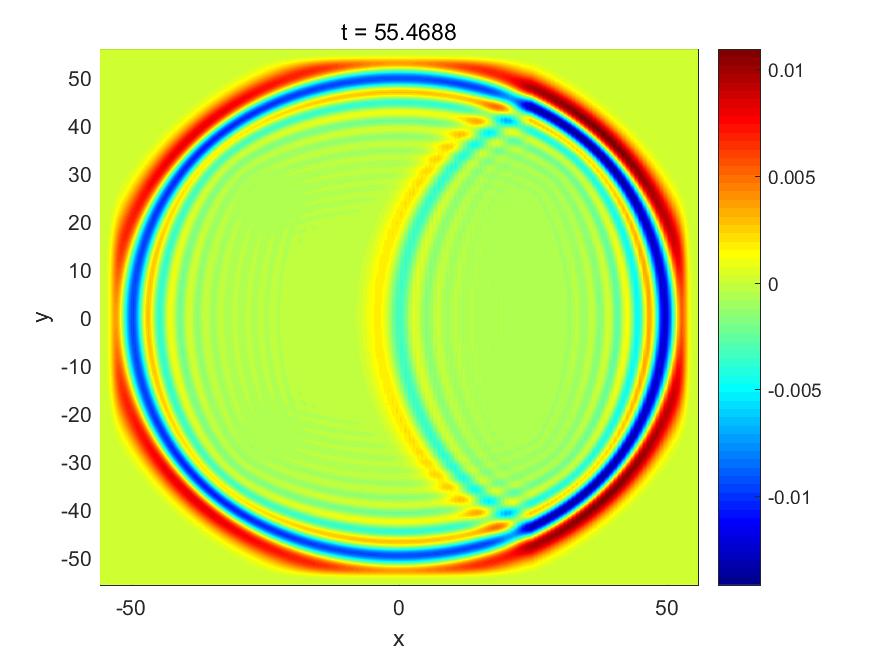}
	\end{minipage} 
    \begin{minipage}{0.28\linewidth}
		\centering
		\includegraphics[width=1\linewidth]{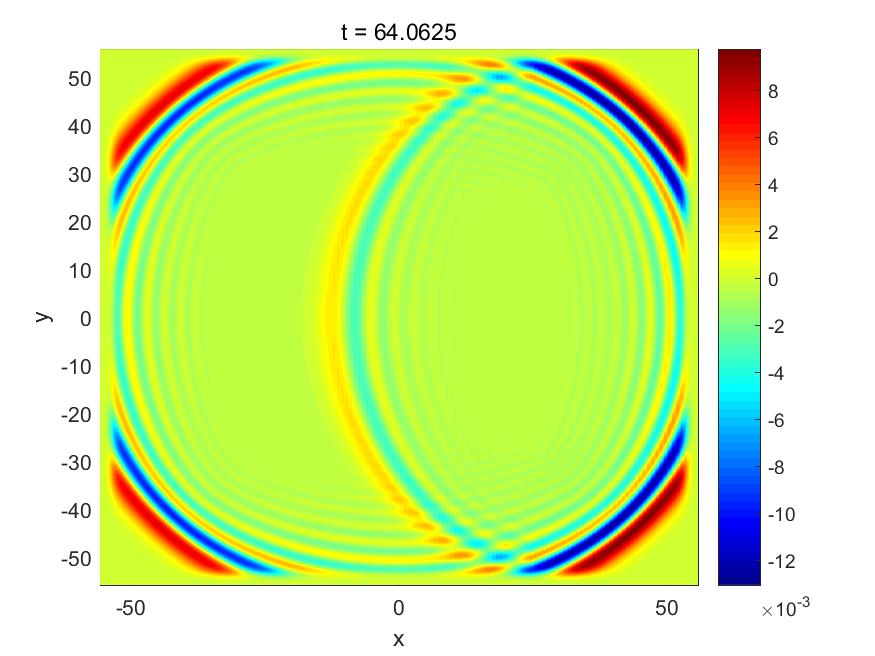}
	\end{minipage}
    \begin{minipage}{0.28\linewidth}
		\centering
		\includegraphics[width=1\linewidth]{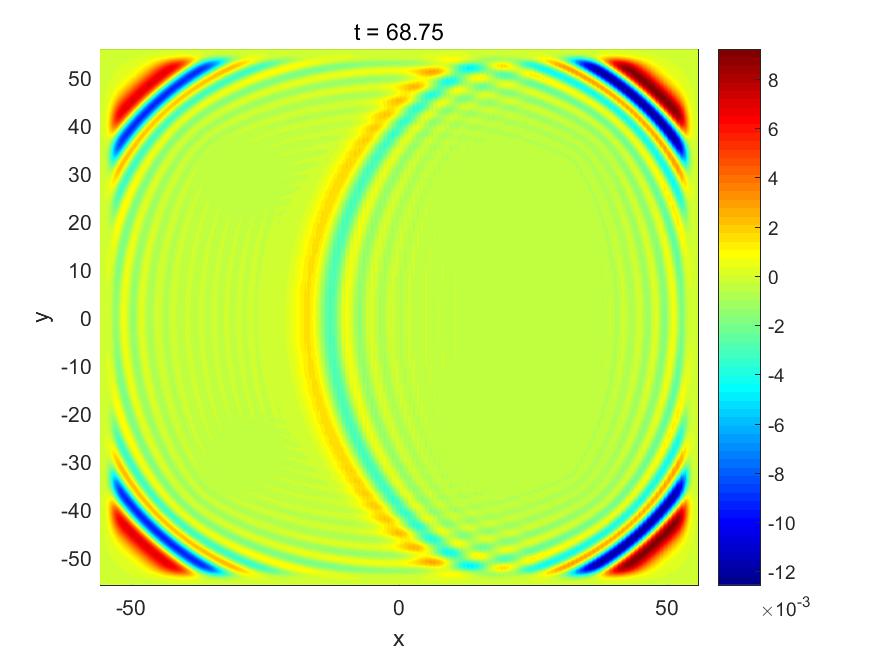}
	\end{minipage} 
 \caption{The snapshot of wave propagation with interface}
 \label{The snapshot of wave propagation with interface}
\end{figure}

The snapshots of the acoustic wave propagation process are shown in Figure \ref{The snapshot of wave propagation with interface}. As can be seen from the figure, the waves are reflected at the interface and are also quickly absorbed in the perfectly matched layer.

\section*{Acknowledgements}
This work is partially supported by the Shanghai Municipal Science and Technology Project 22JC1401600. 
TY was also supported by the Science and Technology Research Project of Jiangxi Provincial Department of Education of China (No. GJJ211027). LZ was also supported by NSFC grant 12271360, and the Fundamental Research Funds for the Central Universities.

\end{document}